# Chains of Integrally Closed Ideals

Kei-ichi Watanabe

ABSTRACT. Let $(A, \mathfrak{m})$ be an excellent normal local ring with algebraically closed residue class field. Given integrally closed $\mathfrak{m}$-primary ideals $I \supset J$, we show that there is a composition series between $I$ and $J$, by integrally closed ideals only. Also we show that any given integrally closed $\mathfrak{m}$-primary ideal $I$, the family of integrally closed ideals $J \subset I$, $l_A(I/J) = 1$ forms an algebraic variety with dimension $\dim A - 1$.

## 1. Introduction

In this paper we investigate chains and familys of integrally closed $\mathfrak{m}$-primary ideals in a Noetherian local ring $(A, \mathfrak{m})$ with algebraically closed residue field.

First question is, "How many integrally closed ideals can we put between two integrally closed $\mathfrak{m}$-primary ideals $I \supset J$" ? The answer is surprisingly simple. In an excellent normal local ring, we can construct a composition series between $I$ and $J$ consisting of integrally closed ideals.

The second question is, "How does the family of integrally closed ideals look like in the family of all ideals of given colength ?"

One surprising fact is that given an integrally closed ideal $I$, the family

$$\{J \subset I \mid l_A(I/J) = 1\}$$

forms an algebraic variety of dimension $\dim A - 1$. On the contrary, the family

$$\{J \supset I \mid l_A(J/I) = 1\}$$

is quite different. For some $I$, that family has arbitrary big dimension while for other $I$, it happens the family consists of one element. The author believes this family reflects some property of the ideal $I$.

When the ring is a rational singularity of dimension 2, the family
$\{J \subset I \mid l_A(I/J) = 1\}$ has a very explicit geometric meaning. We treat this fact in Section 4.

2000 *Mathematics Subject Classification.* Primary 13B22; Secondary 13A30, 14N99.

This work is partially supported by Grants in Aid in sciencific researches, 13440015, 13874006 and also by the project "Study of Singularities", Natural Sciences Research Institute, Nihon University.

I apologize that my proofs are very simple and I don't have nice applications yet. But I believe the phenomena themselves are interesting enough to look at.

## 2. Chains of integrally closed ideals

The main result of this section is the following.

THEOREM 2.1. *Let $(A, \mathfrak{m})$ be an excellent normal local ring with $k = A/\mathfrak{m}$ algebraically closed and let $I \supset J$ be $\mathfrak{m}$-primary integrally closed ideals. Then there exists an integrally closed ideal $I'$ with $I \supset I' \supset J$ and $l_A(I/I') = 1$.*

PROOF. Let $f \colon X \to \mathrm{Spec}(A)$ be a projective birational morphism, where $X$ is normal and $IO_X$ and $JO_X$ are both invertible (for example, $X$ is the normalization of the blowing-up of $IJ$). Then put $IO_X = O_X(-Z)$ and $JO_X = O_X(-Z - Y)$, where $Z, Y$ are effective Cartier divisors on $X$ so that $I = H^0(X, O_X(-Z))$ and $J = H^0(X, O_X(-Z - Y))$. Now, since $A$ is excellent and $X$ is normal, $\mathrm{supp}\,(Y)$ contains a smooth closed point $P$ of $X$. Now, let $g \colon X' \to X$ be the blowing-up of $P$. Since $P$ is a smooth closed point of $X$, $E \colon= g^{-1}(P) \cong \mathbb{P}_k^{d-1}$ ($d = \dim A$). Also, since $P \in \mathrm{supp}\,(Y)$, we have $g^*(Z) < g^*(Z) + E \leq g^*(Z + Y)$.

Now, consider the following exact sequence;

$$0 \to O_{X'}(-g^*(Z) - E) \to O_{X'}(-g^*(Z)) \to O_{X'}(-g^*(Z)) \otimes_{O_{X'}} O_E \to 0.$$

Then we have $O_{X'}(-g^*(Z)) \otimes_{O_{X'}} O_E \cong O_E$ and since $O_X(-g^*(Z))$ is generated by its global sections, we have an exact sequence of global sections

$$0 \to H^0(X, O_X(-g^*(Z) - E)) \to H^0(X, O_X(-g^*(Z))) \to H^0(E, O_E) \to 0.$$

Since $H^0(X, O_X(-g^*(Z))) = I$ and $H^0(E, O_E) \cong k$, $I' \colon = H^0(X, O_X(-g^*(Z) - E))$ is the desired ideal. In fact, $g^*(Z) + E \leq g^*(Z + Y)$ implies $I' \supset J$. □

COROLLARY 2.2. *Let $(A, \mathfrak{m})$ be as in 2.1 and let $I \supset J$ be $\mathfrak{m}$-primary ideals of $A$. Then there is a sequence of $\mathfrak{m}$-primary integrally closed ideals*

$$I = I_0 \supset I_1 \supset \ldots \supset I_l = J$$

*such that $l_A(I_i/I_{i+1}) = 1$ for every $i$.*

## 3. Family of integrally closed ideals

The proof above suggests that there are plenty integrally closed ideals $I'$ with the desired property. In fact, we can show the following.

PROPOSITION 3.1. *Let $(A, \mathfrak{m})$ be any Noetherian local ring with $k = A/\mathfrak{m}$ algebraically closed and $I$ be an integrally closed $\mathfrak{m}$-primary ideal of $A$. Then there is a one-to-one correspondence between the set of integrally closed ideals $I'$ with $I \supset I'$, $l_A(I/I') = 1$ and the set of closed points of $\mathrm{Proj}(\oplus_{n \geq 0}(I^n/\mathfrak{m}I^n))$.*

*In particular, the set of such ideals forms an algebraic variety of dimension $d - 1$.*

PROOF. Put $R \colon = \oplus_{n \geq 0}(I^n/\mathfrak{m}I^n)$. Then $R$ is generated by $V \colon = I/\mathfrak{m}I$ over $R_0 = k$. An ideal $I' \subset I$ with $l_A(I/I') = 1$ is determined by a linear subspace $W \subset V$ of codimension 1 so that $(I' + \mathfrak{m}I)/\mathfrak{m}I = W$. Now, put $R = S/\mathfrak{a}$, where $S = S_k(V)$ be the symmetric algebra over $V$. By definition, the closed points of $\mathrm{Proj}(R)$ is in one to one correspondence with codimension 1 subspaces $W \subset V$ with $WR \supset \mathfrak{a}$. Now let $I' \subset I$ be an ideal generated by $\mathfrak{m}I$ and $W \subset V$, $\dim V/W = 1$. It suffices to show that $I'$ is integrally closed if and only if $WS \supset \mathfrak{a}$.

Now, let $(y_1, \ldots, y_s)$ be a basis of $W$ and $I = (I', x)$. Let $(Y_1, \ldots, Y_s, X)$ be corresponding generators of $S$. If $WS \not\supset \mathfrak{a}$, then let $F = F(X, Y_1, \ldots, Y_s) \in \mathfrak{a}$ be a homogeneous polynomial of degree, say, $n$ not in $WS$. Then we may assume $F$ has the term $X^n$. Since $F \in \mathfrak{a}$, $F(x, y_1, \ldots, y_s) \in \mathfrak{m}I^n$, which shows that $x$ is integral over $I'$.

Conversely, assume $I'$ is not integrally closed. Since $I$ is integrally closed, $I$ is the integral closure of $I'$. Hence $x$ is integral over $I'$ and the integral equation of $x$ over $I'$ gives a homogeneous polynomial in $\mathfrak{a}$ not contained in $(Y_1, \ldots, Y_s)$. □

If $\dim A = d$, $\dim \operatorname{Proj}(\oplus_{n \geq 0}(I^n/\mathfrak{m}I^n)) = d - 1$ for every $\mathfrak{m}$-primary ideal $I$. So, "Looking down from above", the family of "adjacent integrally closed ideals" behaves very uniformly.

At first sight, one might think 3.1 suggest that the set of integrally closed ideals of given colength forms a Zariski closed subset in the family of all ideals. But that is not the case as 3.5 shows.

CONJECTURE 3.2. Let $(A, \mathfrak{m})$ be a Noetherian local ring of dimension $d$ with $k = A/\mathfrak{m}$ algebraically closed. Then, given $n$, the family of integrally closed ideals of colength $n$ is a finite union of locally closed subsets in the family of all the ideals of colength $n$ and a component with maximal dimension has dimesnsion $(d-1)(n-1)$.

But "Looking above from downstairs", the family of "adjacent integrally closed ideals" is very far from uniform as is shown in the following examples.

EXAMPLE 3.3. Let $(A, \mathfrak{m})$ be a local ring of dimension $d$ with $A/\mathfrak{m}$ algebraically closed. Assume that $A$ has a filtration $\mathrm{F} = \{F_i\}_{i \geq 0}$ with $F_1 = \mathfrak{m}, F_i F_j \subset F_{i+j}$ for every $i, j$ and $G_\mathrm{F} = \oplus_{i \geq 0} F_i/F_{i+1}$ is an *integral domain*.

Let $I$ be an ideal such that $F_i \subset I \subset F_{i-1}$ and $l_A(I/F_i) = 1$. Then $I$ is integrally closed. For example, if $A$ is a regular local ring and $I \supset \mathfrak{m}^n$ with $l_A(I/\mathfrak{m}^n) = 1$, then $I$ is integrally closed. In particular, given a $\mathfrak{m}$-primary ideal $\mathfrak{a}$, the family

$$\{J | J \supset \mathfrak{a}, l_A(I/\mathfrak{a}) = 1, \text{ and integrally closed}\}$$

can have arbitrary big dimension.

PROOF. Let $I$ be an ideal with $F_{i-1} \supset I \supset F_i$ with $l_A(I/F_i) = 1$. Put $I = (F_i, x)$ for some $x \in F_{i-1}$. Note that since $G_\mathrm{F}$ is an integral domain, every $F_i$ is integrally closed.

Now, assume $I$ is *not* integrally closed. Recall that the integral closure of $I$ is contained in $F_{i-1}$. Take $y \in F_{i-1}$, which is integral over $I$ and not in $I$. Then we have an integral equation

$$y^n + a_1 y^{n-1} + \ldots + a_{n-1} y + a_n = 0 \quad (a_s \in I^s, 0 \leq s \leq n-1).$$

Note that every term of the equation above is in $F_{n(i-1)}$. Now, considering this equation modulo $F_{n(i-1)+1}$, we have a homogeneous equation

$$Y^n + c_1 X Y^{n-1} + \ldots + c_{n-1} X^{n-1} Y + c_n X^n = 0,$$

where we put $Y = y \mod F_i$, $X = x \mod F_i$ and we put $a_s \equiv c_s x^s \mod F_{s(i-1)+1}$ with some $c_s \in A/\mathfrak{m}$. Since a homogeneous form over an algebraically closed field factors into linear terms and since $G_\mathrm{F}$ is an integral domain, we have $y - dx \in F_i$ for some $d \in A$. This implies $y \in I$, contradicting our assumption. □

EXAMPLE 3.4. Let $(A, \mathfrak{m})$ be a regular local ring of dimension 2 with $A/\mathfrak{m}$ algebraically closed and $I$ be an irreducible (that is, $I$ does not factor as a product $JK$ of integrally closed ideals $J, K$) $\mathfrak{m}$-primary ideal. Then by [4] [5], $I$ has *unique* integrally closed adjacent ideal $J \supset I$. (Conversely, in this case $I$ has unique integrally closed adjacent ideal $J \supset I$ only if $I$ is irreducible ([6])).

EXAMPLE 3.5. Let $(A, \mathfrak{m})$ be a regular local ring with $k = A/\mathfrak{m}$ algebraically closed. Then every $\mathfrak{m}$-primary ideal $I$ of colength $\leq 3$ is integrally closed.

Let $A$ be a regular local ring of dimension 2 and consider the family $X$ of all $\mathfrak{m}$-primary ideals $I$ of colength 4. Then by [1], [2], the family is an irreducible rational variety of dimension 3 and the ideals $I$ not contained in $\mathfrak{m}^2$ form an open dense subvariety $U$. The complement $Y := X \setminus U$ is isomorphic to $\mathbb{P}_k^2$ and corresponds to ideals $I$ contained in $\mathfrak{m}^2$. The family of ideals which are complete intersections forms an open subvariety of $Y$, whose complement (corresponding to integrally closed ideals in $\mathfrak{m}^2$) forms a line in $Y$.

Thus the family of integrally closed ideals of colength 4 is a union of $U$ and a line in $Y$ (closed) and is not Zariski closed in $X$.

## 4. Case of 2-dimensional rational singularities

Let us see the family of integrally closed ideals more explicitly in the case $A$ is a rational singularity of dimension 2. In this section $(A, \mathfrak{m})$ is a 2-dimensional rational singularity with $k = A/\mathfrak{m}$ algebraically closed.

We prepare some terminology on resolution of singularities of $A$.

The following fact is fundamental in our argument.

PROPOSITION 4.1. ([3]) *Let $I$ be an integrally closed $\mathfrak{m}$-primary ideal and take a resolution $f: X \to \mathrm{Spec}(A)$ of $A$ such that $I\mathcal{O}_X$ is invertible and $Z$ be a divisor on $X$ with $I\mathcal{O}_X \cong \mathcal{O}_X(-Z)$. Then for every irreducible exceptional curve $E$ on $X$, $EZ \leq 0$. We call an exceptional divisor $Z$ with this property an anti-nef divisor.*

*Conversely, if $Z$ is an anti-nef divisor, then $I := H^0(X, \mathcal{O}_X(-Z))$ is an integrally closed ideal and $\mathcal{O}_X(-Z) = I\mathcal{O}_X$. Thus there is a one-to-one correspondence between the set of anti-nef cycles on $X$ and integrally closed $\mathfrak{m}$-primary ideals such that $I\mathcal{O}_X$ is invertible.*

The following statement gives a geometric description of integrally closed ideals $J$ with $l_A(I/J) = 1$ for a fixed integrally closed ideal $I = H^0(X, \mathcal{O}_X(-Z))$.

PROPOSITION 4.2. *Let $A$ be a 2-dimensional rational singularity and let $f: X \to \mathrm{Spec}(A)$ be as in previous proposition and $I = H^0(X, \mathcal{O}_X(-Z))$, where $Z$ is an anti-nef divisor on $X$. If $J \subset I$ is an integrally closed ideal such that $l_A(I/J) = 1$, then either*

(1) $J = H^0(X, \mathcal{O}_X(-Z-Y))$, *where $Y$ is an effective exceptional cycle on $X$ such that $Z + Y$ is anti-nef satisfying the following conditions;*
   (a) $ZY = 0$ *and* (b) $p_a(Y) = 0$,
   *where $p_a(Y)$ is defined by* $p_a(Y) = \dfrac{Y^2 + K_X Y}{2} + 1$.

(2) $J = H^0(X', \mathcal{O}_{X'}(-\pi^*(Z) - E))$, *where $\pi: X' \to X$ is a blowing-up of a point $P \in f^{-1}(\mathfrak{m})$ such that every irreducible exceptional curve $C \ni P$ satisfies $CZ < 0$.*

PROOF. If $I = H^0(X, \mathcal{O}_X(-Z))$, where $Z$ is an anti-nef exceptional divisor on $A$, we have the following formula by the "Theorem of Riemann-Roch" on $X$ (cf. [3]);

$$l_A(A/I) = -\frac{Z^2 + K_X \cdot Z}{2}.$$

Hence if $J = H^0(X, \mathcal{O}_X(-Z-Y))$ with $Z+Y$ anti-nef, then $l_A(I/J) = -\frac{(Z+Y)^2 + K_X \cdot (Z+Y)}{2} + \frac{Z^2 + K_X \cdot Z}{2} = -Z \cdot Y + (1 - p_a(Y))$.

Now, $Z \cdot Y \leq 0$, since $Z$ is anti-nef and since $A$ is a rational singularity, $p_a(Y) \leq 0$ for every exceptional cycle on $X$. So, if $l_A(I/J) = 1$, we must have the condition (1).

Now, an integrally closed ideal $J \subset I$ is written as $J = H^0(X', \mathcal{O}_{X'}(-g^*(Z) - Y))$, where $g \colon X' \to X$ be composite of blowing-ups of points and $g^*(Z) + Y$ is an anti-nef divisor on $X'$. If $g(Y)$ is a divisor on $X$, then $Z + g_*(Y)$ is also anti-nef being the image of an anti-nef divisor and since $J = H^0(X', \mathcal{O}_{X'}(-g^*(Z) - Y)) \subset H^0(X, \mathcal{O}_X(-Z - g_*(Y))) \subset I$, $J = H^0(X, \mathcal{O}_X(-Z - g_*(Y)))$ and we must be in case (1).

Now, assume $\dim g(Y) = 0$. By (1), we have $Y \cdot g^*(Z) = 0$ and $p_a(Y) = 0$. By the formula $p_a(Y + Y') = p_a(Y) + p_a(Y') + Y \cdot Y' - 1$, the condition $p_a(Y) = 0$ imposes $Y$ is connected. Hence $g(Y) = P$ is a point on $X$. Let $\pi_1 \colon X_1 \to X$ be the blowing-up of the point $P$ and factor $g \colon X' \to X$ into

$$g \colon X' = X_n \to X_{n-1} \to \ldots \to X_1 \to X,$$

where we denote $\pi \colon X_i \to X_{i-1}$ and $h_i \colon X' \to X_i$. Since $Y$ is connected, we may assume the center of $\pi_i$ is on some exceptional curve of $\pi_j (j < i)$.

Now, if $h_1(Y)$ is again a point, let $i > 1$ be the least number so that $h_i(Y)$ is a divisor on $X_i$ and assume the center of $\pi_i$ lies on the strict transform $E_j$ of the center of $\pi_j, j < i$. Then we have $(h_i)_*(Y) \cdot E_j > 0$ and if $E'$ is the strict transform of $E_j$, we still have $(g^*(Z) + Y) \cdot E' > 0$ or if there is another blowing-up at the intersection of $Y$ and $E_j$ so that $Y$ does not contain the new exceptional curve, then $(g^*(Z) + Y) \cdot E'' > 0$ for a new exceptional curve $E''$ (note that $g^*(Z) \cdot E = 0$ for every strict transform of the exceptional curves of $\pi_i, (i = 1, \ldots, n)$). Anyway, this contradicts $g^*(Z) + Y$ is anti-nef. So, $h_1(Y)$ must contain the exceptional curve $E_1$ of $\pi_1$. Then $h_1^*(E_1) \leq Y$ and $J = H^0(X, \mathcal{O}_X(-g^*(Z) - Y)) \subset J' := H^0(X_1, \mathcal{O}_X(-\pi_1^*(Z) - E_1)) \subset I$. Since $l_A(I/J') = 1$, $J = J'$ and we are in case (2).

Now, if the center of $\pi_1$ is on some exceptional curve $C$ with $Z \cdot C = 0$, then the strict transform $C'$ has $(\pi_1^*(Z) + E_1) \cdot C' = (\pi_1^*(Z) + E_1) \cdot (\pi_1^*(C) - E_1) = 1 > 0$, contradicting the fact $\pi_1^*(Z) + E_1$ should be anti-nef. □

REMARK 4.3. Let $f \colon X \to \mathrm{Spec}(A)$ be as above and $E = \bigcup C_i = f^{-1}(\mathfrak{m})$. Let $E_0 \subset E$ be the union of irreducible component $C_i$ with $Z \cdot C_i = 0$ and $B$ be a connected component of $E_0$. If $Y$ is as in 4.2 (1) and if $Y$ contains an irreducible curve $C \subset B$, then since $Z + Y$ is anti-nef, $Y$ must contain every irreducible curve in $B$. Thus there is a one to one correspondence between cycles $Y$ as in 4.2 (1) and the connected components of $E_0$.

EXAMPLE 4.4. Let $A = k[[X, Y, Z]]/(XY - Z^{n+1})$ and $I = \mathfrak{m}$, the maximal ideal of $A$. Then we can take $X$ to be the minimal resolution $f \colon X \to \mathrm{Spec}(A)$ of $A$, the exceptional curve $E = f^{-1}(\mathfrak{m})$ consists of the chain of $n$ $\mathbb{P}^1$'s with self-intersection number $-2$ and $\mathfrak{m}\mathcal{O}_X = \mathcal{O}_X(-E)$. In this case, $E_0$ is the chain of $n - 2$ middle curves and an integrally closed ideal $J$ of colength 2 is given either $J = H^0(X, \mathcal{O}_X(-E - E_0))$ or by blowing-up a point on one of two end curves.

On the other hand, since $G_{\mathfrak{m}}(A) = \oplus_{n \geq 0} \mathfrak{m}^n/\mathfrak{m}^{n+1} \cong k[X,Y,Z]/(XY)$ $(n \geq 2)$ or $k[X,Y,Z]/(XY - Z^2)$ $(n = 1)$, if $n \geq 2$, $\mathrm{Proj}(G_{\mathfrak{m}}(A))$ consists of two $\mathbb{P}^1$'s intersecting at one point. Thus the points on each $\mathbb{P}^1$ except the intersection corresponds to the two end curves and the intersection point corresponds to the ideal $J = H^0(X, \mathcal{O}_X(-E - E_0))$.

**Acknowledgement.** The author thanks D.Eisenbud and J.Lipman for valuable comments.

## References

[1] J. Briançon, Description de Hilb$^n$ $\mathbb{C}\{x,y\}$, *Inventiones math.*, **41** (1977) 45-89.
[2] A. Iarrobino, Punctual Hilbert schemes, *Bull. A.M.S.*, **78** (1972) 819-823.
[3] J. Lipman, Rational singularities with applications to algebraic surfaces and unique factorization., *Publ. I.H.E.S*, **36** (1969) 195–279.
[4] J. Lipman, Proximity inequalities for complete ideals in two-dimensional regular local rings, *Contemporary Math.*, **159** (1994) 293-306.
[5] S. Noh, Adjacent integrally closed ideals in dimension two, *J. P. A. A.*, **85** (1993), 163-184.
[6] S. Noh and K.-i. Watanabe, Adjacent integrally closed ideals in two dimensional regular local rings, in preparation.

DEPARTMENT OF MATHEMATICS, COLLEGE OF HUMANITY AND SCIENCES, NIHON UNIVERSITY, SETAGAYA-KU, TOKYO 156–8550, JAPAN

*E-mail address*: `watanabe@math.chs.nihon-u.ac.jp`